\documentclass{amsart}

\usepackage{amssymb}

\newtheorem{lemma}{Lemma}[section]
\newtheorem{theorem}[lemma]{Theorem}
\newtheorem{proposition}[lemma]{Proposition}
\newtheorem{corollary}[lemma]{Corollary}
\newtheorem{problem}[lemma]{Problem}
\newtheorem{question}[lemma]{Question}
\theoremstyle{definition}
\newtheorem{definition}[lemma]{Definition}
\newtheorem{example}[lemma]{Example}
\newtheorem{remark}[lemma]{Remark}

\newcommand{\F}{\mathcal {F}}
\newcommand{\Ra}{\Rightarrow}
\newcommand{\U}{\mathcal U}
\newcommand{\V}{\mathcal V}
\newcommand{\A}{\mathcal A}
\newcommand{\w}{\omega}
\newcommand{\IN}{\mathbb N}
\newcommand{\IZ}{\mathbb Z}
\newcommand{\C}{\mathcal C}

\newcommand{\LL}{\mathcal L}

\title{Algebra in superextension of groups, II: cancelativity and
 centers}
\author{Taras Banakh and Volodymyr Gavrylkiv}
\address{Ivan Franko National University of Lviv, Ukraine
 and\newline\phantom{ }\hskip8pt Akademia \'Swi\c etokrzyska, Kielce, Poland}
\email{tbanakh@yahoo.com}
\address{Vasyl Stefanyk Precarpathian National University,
 Ivano-Frankivsk, Ukraine}
\email{vgavrylkiv@yahoo.com} \subjclass{20M99, 54B20}

\begin{document}
\begin{abstract}
Given a countable group $X$ we study the algebraic structure of its
 superextension
 $\lambda(X)$. This is a right-topological semigroup consisting of all
 maximal linked systems
 on $X$ endowed with the operation $$\A\circ\mathcal B=\{C\subset
 X:\{x\in X:x^{-1}C\in\mathcal B\}\in\A\}$$
 that extends the group operation of $X$. We show that the subsemigroup
 $\lambda^\circ(X)$ of free maximal linked systems contains an open
 dense subset of right cancelable elements. Also we prove that the
 topological center of $\lambda(X)$ coincides with the subsemigroup
 $\lambda^\bullet(X)$ of all maximal linked systems with finite support. This
 result is applied to show that the algebraic center of $\lambda(X)$
 coincides with the algebraic center of $X$ provided $X$ is countably infinite.
 On the other hand, for finite groups $X$ of order $3\le|X|\le5$ the
 algebraic center of $\lambda(X)$ is strictly larger than the algebraic
 center of $X$.
\end{abstract}
\maketitle

\section*{Introduction}

 After the topological proof (see \cite[p.102]{HS}, \cite{H2}) of Hindman theorem \cite{Hind}, topological methods become a standard tool in the modern
combinatorics of numbers, see \cite{HS}, \cite{P}. The crucial point is that any semigroup operation $\ast$ defined on
any discrete space $X$ can be extended to a right-topological semigroup operation on $\beta(X)$, the Stone-\v Cech
compactification of $X$. The extension of the operation from $X$ to $\beta(X)$ can be defined by the simple formula:
\begin{equation}\label{extension}\U\ast\V=\Big\{\bigcup_{x\in
 U}x*V_x:U\in\U,\;\{V_x\}_{x\in U}\subset\V\Big\},
\end{equation}
where $\U,\V$ are ultrafilters on $X$. Endowed with the so-extended operation, the Stone-\v Cech compactification
$\beta(X)$ becomes a compact right-topological semigroup. The algebraic properties of this semigroup (for example, the
existence of idempotents or minimal left ideals) have important consequences in combinatorics of numbers, see
\cite{HS}, \cite{P}.

The Stone-\v Cech compactification $\beta(X)$ of $X$ is the subspace of the double power-set $\mathcal P(\mathcal
P(X))$, which is a complete lattice with respect to the operations of union and intersection. In \cite{G2} it was
observed that the semigroup operation extends not only to $\beta(X)$ but also to the complete sublattice $G(X)$ of
$\mathcal P(\mathcal P(X))$ generated by $\beta(X)$. This complete sublattice consists of all inclusion hyperspaces
over $X$.

By definition, a family $\F$ of non-empty subsets of a discrete space $X$ is called an {\em inclusion hyperspace} if
$\F$ is monotone in the sense that a subset $A\subset X$ belongs to $\F$ provided $A$ contains some set $B\in\F$. On
the set $G(X)$ there is an important transversality operation assigning to each inclusion hyperspace $\F\in G(X)$ the
inclusion hyperspace
$$\F^\perp=\{A\subset X:\forall F\in\F\;(A\cap F\ne\emptyset)\}.$$
This operation is involutive in the sense that $(\F^\perp)^\perp=\F$.

It is known that the family $G(X)$ of inclusion hyperspaces on $X$ is closed in the double power-set $\mathcal
P(\mathcal P(X))=\{0,1\}^{\mathcal P(X)}$ endowed with the natural product topology.

The extension of a binary operation $\ast$ from $X$ to $G(X)$ can be defined in the same way as for ultrafilters, i.e.,
by the formula~(\ref{extension}) applied to any two inclusion hyperspaces $\U,\V\in G(X)$. In \cite{G2} it was shown
that for an associative binary operation $\ast$ on $X$ the space $G(X)$ endowed with the extended operation becomes a
compact right-topological semigroup. Besides the Stone-\v Cech extension, the semigroup $G(X)$ contains many important
spaces as closed subsemigroups. In particular,  the space
$$\lambda(X)=\{\F\in G(X):\F=\F^\perp\}$$of maximal linked systems on
 $X$ is a closed subsemigroup of $G(X)$. The space $\lambda(X)$ is
  well-known in  General and Categorial Topology as the {\em
 superextension}
  of $X$, see \cite{vM}, \cite{TZ}.  Endowed with the extended binary
 operation, the
  superextension $\lambda(X)$ of a semigroup $X$ is a supercompact
  right-topological semigroup containing $\beta(X)$ as a subsemigroup.

The thorough study of algebraic properties of the superextensions
of
 groups was started in \cite{BGN}
where we described right and left zeros in $\lambda(X)$ and detected
 all groups $X$ with commutative superextension $\lambda(X)$ (those are
 groups of cardinality $|X|\le 4$). In \cite{BGN} we also described the
 structure of the semigroups $\lambda(X)$ for all finite groups $X$ of
 cardinality $|X|\le 5$.
In \cite{BG3} we shall describe the structure of minimal left ideals of
 the superextensions of groups. In this paper we  concentrate at
  cancellativity and centers (topological and algebraic) in the
 superextensions $\lambda(X)$ of groups $X$. Since $\lambda(X)$ is an intermediate
 subsemigroup between $\beta(X)$ and $G(X)$ the obtained results for
 $\lambda(X)$ in a sense are intermediate between those for $\beta (X)$ and
 $G(X)$.

In section 2 we describe cancelable elements of $\lambda(X)$. In
 particular, we show that for a finite group $X$ all left or right cancelable
 elements of $\lambda(X)$ are principal ultrafilters. On the other hand,
 if a group $X$ is countable, then the set of right cancelable elements
 has open dense intersection with the  subsemigroup
 $\lambda^\circ(X)\subset\lambda(X)$ of free maximal linked systems, see
 Theorem~\ref{t2.4}. This resembles the situation with the semigroup $\beta(X)\setminus X$ which
 contains a dense open subset of right cancelable elements (see
 \cite[8.10]{HS}), and also with the semigroup $G(X)$ whose right cancelable elements
 form a subset having open dense intersection with the set $G^\circ(X)$
 of free inclusion hyperspaces, see \cite{G2}.

The section 3 is devoted to describing the topological center of
 $\lambda(X)$.
By definition, the {\em topological center} of a right-topological semigroup $S$ is the set $\Lambda(S)$ of all
elements $a\in S$ such that the left shift $l_a:S\to S$, $l_a:x\mapsto a*x$, is continuous. By \cite{HS} for every
group $X$ the topological center of the semigroup $\beta(X)$ coincides with $X$. On the other hand, the topological
center of the semigroup $G(X)$ coincides with the subspace $G^\bullet(X)$ of $G(X)$ consisting of inclusion hyperspaces
with finite support, see \cite[7.1]{G2}. A
 similar results holds also for the semigroup $\lambda(X)$: for any
at most countable group $X$ the topological center of $\lambda(X)$ coincides with $\lambda^\bullet(X)$, see
Theorem~\ref{t3.4}.

The final section \ref{center} is devoted to describing the algebraic center of $\lambda(X)$. We recall that the {\em
algebraic center} of a semigroup $S$ consists of all elements $s\in S$ that commute with all other elements of $S$. In
Theorem~\ref{algcent} we shall prove
 that
for any countable infinite group $X$ the algebraic center of $\lambda(X)$ coincides with the algebraic center of $X$.
It is interesting to note that for any group $X$ the algebraic centers of the semigroups
 $\beta(X)$
and $G(X)$ also coincide with the center of the group $X$, see
 \cite[6.54]{HS} and  \cite[6.2]{G2}. In contrast, for finite groups
 $X$ of cardinality $3\le|X|\le5$
the algebraic center of $\lambda(X)$ is strictly larger than the
 algebraic center of $X$, see Remark~\ref{r4.4}.

\section{Inclusion hyperspaces and superextensions}

In this section we recall the necessary definitions and facts.

A family $\mathcal L$ of subsets of a set $X$ is called {\em a linked
 system on $X$}
 if $A\cap B\ne\emptyset$ for all $A,B\in\mathcal L$ and $\mathcal L$
 is closed under taking supersets.
 Such a linked system $\mathcal L$
 is {\em maximal linked} if $\mathcal L$ coincides with any linked
 system $\mathcal L'$
 on $X$ that contains $\mathcal L$. Each (ultra)filter on $X$ is a
 (maximal) linked system.
By $\lambda(X)$ we denote the family of all maximal linked systems on $X$. Since each ultrafilter on $X$ is a maximal
linked system, $\lambda(X)$ contains the  Stone-\v Cech extension $\beta(X)$ of $X$. It is easy to see that each
maximal linked system on $X$ is an inclusion hyperspace on $X$ and hence $\lambda(X)\subset G(X)$. Moreover, it can be
shown that $\lambda(X)=\{\A\in G(X):\A=\A^\perp\}$, see \cite{G1}.

By \cite{G1} the subspace $\lambda(X)$ is closed in the space $G(X)$ endowed with the topology generated by the
sub-base consisting of the sets
$$U^+=\{\A\in G(X):U\in\A\}\mbox{ and }U^-=\{\A\in
 G(X):U\in\A^\perp\}$$
where $U$ runs over subsets of $X$. By \cite{G1} and \cite{vM} the spaces $G(X)$ and $\lambda(X)$ are supercompact in
the sense that any their cover by the sub-basic sets contains a two-element subcover. Observe that
$U^+\cap\lambda(X)=U^-\cap\lambda(X)$ and hence the topology on $\lambda(X)$ is generated by the sub-basis consisting
of the sets
$$U^\pm=\{\A\in \lambda(X):U\in\A\}, \;\;U\subset X.$$

We say that an inclusion hyperspace $\A\in G(X)$
\begin{itemize}
\item has {\em finite support} if there is a finite family
 $\F\subset\A$ of
 finite subsets of $X$ such that each set $A\in\A$ contains a set
 $F\in\F$;
\item is {\em free} if for each $A\in\A$ and each finite subset
$F\subset X$ the complement $A\setminus F$ belongs to $\A$.
\end{itemize}

By $G^\bullet(X)$ we denote the subspace of $G(X)$ consisting of inclusion hyperspaces with finite support and
$G^\circ(X)$ stands for the subset of free inclusion hyperspaces on $X$. Those two sets induce the subsets
$$\lambda^\bullet(X)=G^\bullet(X)\cap\lambda(X)\mbox{ and
 }\lambda^\circ(X)=G^\circ(X)\cap \lambda(X)$$
in the superextension $\lambda(X)$ of $X$. By \cite{G1}, $\lambda^\bullet(X)$ is an open dense subset of $\lambda(X)$
while $\lambda^\circ(X)$ is  closed and nowhere dense in $\lambda(X)$.

Given any semigroup operation $\ast:X\times X\to X$ on a set $X$ we can
 extend this operation to $G(X)$ letting
$$\U\ast\V=\Big\{\bigcup_{x\in U}x*V_x:U\in\U,\;\{V_x\}_{x\in
 U}\subset\V\Big\}$$for inclusion
hyperspaces $\U,\V\in G(X)$. Equivalently, the product $\U\ast\V$ can be defined as
\begin{equation}\label{eq2}
\U\ast\V=\{A\subset X:\{x\in X:x^{-1}A\in \V\}\in\U\}
\end{equation}
where $x^{-1}A=\{z\in X:x\ast z\in A\}$. By \cite{G2} the so-extended operation turns $G(X)$ into a right-topological
semigroup. The structure of this semigroup was studied in details in \cite{G2}. In this paper we shall concentrate at
the study of the algebraic structure of the semigroup $\lambda(X)$ for a group $X$.

The formula (\ref{eq2}) implies that the product $\U\ast\V$ of two maximal linked systems $\U$ and $\V$ is a principal
ultrafilter if and only if both $\U$ and $\V$ are principal ultrafilters. So we get the following

\begin{proposition}\label{ideal} For any group $X$ the set
 $\lambda(X)\setminus X$ is a two-sided ideal in $\lambda(X)$.
\end{proposition}

\section{Cancelable elements of $\lambda(X)$}\label{cancelable}

 In this section, given a group $X$ we shall detect  cancelable
elements of $\lambda(X)$.

We recall that an element $x$ of a semigroup $S$ is {\em right} (resp.
 {\em left}) {\em cancelable}
 if for every $a,b\in X$ the equation $x*a=b$ (resp. $a*x=b$) has at
 most one solution $x\in S$.
This is equivalent to saying that the right (resp. left) shift
 $r_a:S\to S$, $r_a:x\mapsto x*a$,
(resp. $l_a:S\to S$, $l_a:x\mapsto a*x$) is injective.

\begin{proposition}\label{p2.1} Let $G$ be a finite group. If
 $\C\in\lambda(G)$ is left or right cancelable,
then $\C$ is a principal ultrafilter.
\end{proposition}

\begin{proof} Assume that some maximal linked system
 $a\in\lambda(G)\setminus G$ is left cancelable.
This means that the left shift $l_a:\lambda(G)\to\lambda(G)$, $l_a:x\mapsto a\circ x$, is injective. By
Proposition~\ref{ideal}, the set $\lambda(G)\setminus G$ is an ideal in $\lambda(G)$. Consequently,
$l_a(\lambda(G))=a*\lambda(G)\subset\lambda(G)\setminus G$. Since $\lambda(G)$ is finite, $l_a$ cannot be injective.
\end{proof}

Thus the semigroups $\lambda(X)$ can have non-trivial cancelable
elements only for infinite groups $X$. According to
\cite[8.11]{HS} an ultrafilter $\U\in\beta(X)$ is right cancelable
if and only if the orbit $\{x\U:x\in X\}$ is discrete in
$\beta(X)$ if and only if for every $x\in X$ there is a set
$U_x\in\U$ such that the indexed family $\{x*U_x:x\in X\}$ is
disjoint.

This characterization admits a partial generalization to the semigroup $G(X)$. According to \cite{G2} if an inclusion
hyperspace $\A\in G(X)$ is right cancelable in $G(X)$, then its orbit $\{x*\A:x\in X\}$ is discrete in $G(X)$. On the
other hand, $\A$ is cancelable provided for every $x\in X$ there is a set $A_x\in\A\cap\A^\perp$ such that the indexed
family $\{x*A_x:x\in X\}$ is disjoint. The latter means that $x*A_x\cap y*A_y=\emptyset$ for any distinct points
$x,y\in X$. This result on right cancelable elements in $G(X)$ will help us to prove a similar result on the right
cancelable elements in the semigroup $\lambda(X)$.

\begin{theorem}\label{cancel} Let $X$ be a group and $\LL\in\lambda(X)$
 be a maximal linked system on $X$.
\begin{enumerate}
\item If $\LL$ is right cancelable in $\lambda(X)$, then the orbit
$\{x\LL:x\in X\}$ is discrete in $\lambda(X)$ and $x\LL\ne y\LL$ for all $x,y\in X$. \item $\LL$ is right cancelable in
$\lambda(X)$ provided for every $x\in X$ there is a set $S_x\in\mathcal L$ such that the family $\{x*S_x:x\in X\}$ is
disjoint.
\end{enumerate}
\end{theorem}

\begin{proof} 1. First note that the right cancelativity of a maximal
 linked system $\LL\in \lambda(X)$
 is equivalent to the injectivity of the map
$\mu_X\circ \lambda\bar R_\LL:\lambda(X)\to \lambda(X)$, see \cite{G2}.
 We recall that $\mu_X:\lambda^2(X)\to \lambda(X)$ is the
 multiplication of the monad $\mathbb
\lambda=(\lambda,\mu,\eta)$ while $\bar R_\LL:\beta(X)\to \lambda(X)$ is the Stone-\v Cech extension of the right shift
$R_\LL:X\to \lambda(X)$, $R_\LL:x\mapsto x*\LL$. The map $\bar R_\LL$ certainly is not injective if $R_\LL$ is not an
embedding, which is equivalent to the discreteness of the indexed set $\{x*\LL:x\in X\}$ in $\lambda(X)$.
\smallskip

2. Assume that $\{S_x\}_{x\in X}\subset\LL$ is a family such that
 $\{x*S_x:x\in X\}$ is disjoint.
To prove that $\LL$ is right cancelable, take two maximal linked
 systems $\A,\mathcal B\in \lambda(X)$ with
$\A\circ\LL=\mathcal B\circ \LL$. It is sufficient to show that
 $\A\subset\mathcal B$. Take any set $A\in\A$ and observe that the set
   $\bigcup_{a\in A}aS_a$ belongs to $\A\circ\LL=\mathcal B\circ\LL$.
 Consequently, there is
   a set $B\in\mathcal B$ and a family of sets $\{L_b\}_{b\in
 B}\subset\LL$ such that
   $$\bigcup_{b\in B}bL_b\subset \bigcup_{a\in A}aS_a.$$ It follows
 from $S_b\in\LL$
   that $L_b\cap S_b$ is not empty for every $b\in B$.

Since the sets $aS_a$ and $bS_b$ are disjoint for different $a,b\in X$, the inclusion $$\bigcup_{b\in B}b(L_b\cap
S_b)\subset \bigcup_{b\in B}bL_b\subset \bigcup_{a\in A}aS_a$$ implies $B\subset A$ and hence $A\in\mathcal B$.
\end{proof}

It is interesting to remark that  the first item gives a necessary but
 not sufficient
 condition of the right cancelability in $\lambda(X)$ (in contrast to
 the situation in $\beta(X)$).

\begin{example} By \cite[6.3]{BGN}, the superextension $\lambda(C_4)$ of
 the 4-element cyclic group $C_4$ is isomorphic to the direct product
$C_4\times C_2^1$, where $C_2^1=C_2\cup\{e\}$ is the 2-element cyclic
 group with attached external unit $e$ (the latter means that $ex=xe=x$
 for all $x\in C_2^1$). Consequently, each element of the ideal
 $\lambda(C_4)\setminus C_4$ is not cancelative but has the discrete 4-element
 orbit $\{x\mathcal L:x\in C_4\}$.  In fact all the (left or right)
 cancelable elements of $\lambda(C_4)$ are principal ultrafilters, see
 Proposition~\ref{p2.1}.
 \end{example}

According to \cite[8.10]{HS}, for each infinite group the
semigroup
 $\beta(X)$ contains many right
 cancelable elements. In fact, the set of right cancelable elements
 contains an open dense
 subset of $\beta(X)\setminus X$. A similar result holds also for the
 semigroup $G(X)$ over
  a countable group $X$: the set of right cancelable elements of $G(X)$
 contains an open dense
  subset of the subsemigroup $G^\circ(X)$. Theorem~\ref{cancel} will
 help us to prove a similar
  result for the semigroup $\lambda(X)$.

\begin{theorem}\label{t2.4} For each counatable group $X$ the
 subsemigroup $\lambda^\circ(X)$ of free maximal
linked systems contains an open dense subset consisting of right
 cancelable elements in the semigroup $\lambda(X)$.
\end{theorem}

\begin{proof}
 Let $X=\{x_n:n\in\w\}$ be an injective enumeration of the countable
 group $X$.
 Given a free maximal linked system $\LL\in \lambda^\circ(X)$ and a
 neighborhood $O(\LL)$ of
  $\LL$ in $\lambda^\circ(X)$, we should find a non-empty open subset
 of right cancelable elements in $O(\LL)$.
  Without loss of generality, the neighborhood $O(\LL)$ is of basic
 form:
$$O(\LL)=\lambda^\circ(X)\cap U_0^\pm\cap\dots\cap U_{n-1}^\pm$$
for some sets $U_1,\dots,U_{n-1}$ of $X$. Those sets are infinite because $\LL$ is free. We are going to construct an
infinite set $C=\{c_n:n\in\w\}\subset X$ that has infinite intersection with the sets $U_i$, $i<n$, and such that for
any distinct $x,y\in X$ the intersection $xC\cap yC$ is finite. The points $c_k$, $k\in\w$, composing the set $C$ will
be chosen by induction to satisfy the following conditions:
\begin{itemize}
\item $c_k\in U_j$ where $j=k \mod n$; \item $c_k$ does not belong
to the finite set $$F_k=\{z\in X:\exists i,j\le k\;\exists l<k\;\;(x_iz=x_jc_l)\}.$$
\end{itemize}
It is clear that the so-constructed set $C=\{c_k:k\in\w\}$ has infinite intersection with each set $U_i$, $i<n$. The
choice of the
 points $c_k$ for
$k>j$ implies that $x_iC\cap x_jC\subset \{x_i c_m:m\le j\})$ is
 finite.

Now let $\C$ be a free maximal linked system on $X$ enlarging the linked system generated by the sets $C$ and
$U_0,\dots,U_{n-1}$. It is clear that $\C\in O(\LL)$. Consider the open neighborhood
$$O(\C)=O(\LL)\cap C^\pm$$ of $\C$ in $\lambda^\circ(X)$.

We claim that each maximal linked system $\A\in O(\C)$ is right cancelable in $\lambda(X)$. This will follow from
Proposition~\ref{cancel} as soon as we construct a family of sets $\{A_i\}_{i\in\w}\in \A$ such that $x_iA_i\cap
x_jA_j=\emptyset$ for any numbers $i<j$. Observe that the sets
$$A_i=C\setminus\{x_i^{-1}x_kc_m:k<i,\;m\le i\},\;\;i\in\w,$$
have the required property.
\end{proof}

By \cite[8.34]{HS}, the semigroup $\beta(\IZ)$ contains many free
 ultrafilters that are left cancelable in  $\beta(\IZ)$.
 On the other hand, by \cite[8.1]{G2}, the only left cancelable
 elements of the semigroup $G(\IZ)$ are
  principal ultrafilters.

\begin{problem} Is there a maximal linked system
 $\U\in\lambda(\IZ)\setminus\IZ$ which is left cancelable
in the semigroup $\lambda(\IZ)$?
\end{problem}

\section{The topological center of $\lambda(X)$}

In this section we describe the topological center of the superextension $\lambda(X)$  of a group $X$. By definition,
the {\em topological center} of a right-topological semigroup $S$ is the set $\Lambda(S)$ of all elements $a\in S$ such
that the left shift $l_a:S\to S$, $l_a:x\mapsto a*x$, is continuous.

By \cite[4.24,\;6.54]{HS} for every group $X$ the topological
center of the semigroup $\beta(X)$ coincides with $X$. On the
other hand, the topological center of the
 semigroup
$G(X)$ coincides with $G^\bullet(X)$, see \cite[7.1]{G2}. A similar
 results holds
also for the semigroup $\lambda(X)$: the topological center of
 $\lambda(X)$ coincides with $\lambda^\bullet(X)$ (at least for countable groups $X$).

To prove this result we shall use so-called detecting ultrafilters.

\begin{definition}\label{ultradetect} A free ultrafilter $\mathcal D$
 on a group $X$ is called
{\em detecting} if there is an indexed family of sets $\{D_x:x\in
 X\}\subset \mathcal D$ such that
 for any $A\subset X$
\begin{enumerate}
\item the set $U_A=\bigcup_{x\in A}xD_x$ has the property: $U_A\cup
 yU_A\ne X$ for all $y\in X$;
\item for every $D\in\mathcal D$ the set $\{x\in X:xD\subset U_A\}$ is
 finite and lies in $A$.
\end{enumerate}
\end{definition}

\begin{lemma}\label{l3.2} On each countable group $X$ there is a
 detecting ultrafilter.
\end{lemma}

\begin{proof} Let $X=\{x_n:n\in\w\}$ be an injective enumeration of the
 group $X$ such
that $x_0$ is the neutral element of $X$. For every $n\in\w$ let
 $F_n=\{x_i,x_i^{-1}:i\le n\}$.
  Let $a_0=x_0$ and inductively, for every $n\in\w$ choose an element
 $a_n\in X$ so that
$$a_n\notin F_n^{-1}F_n A_{<n}\mbox{ \ where \ }A_{<n}=\{a_i:i<n\}.$$
For every $n\in\w$ let $A_{\ge n}=\{a_{i}:i\ge n\}$. Let also
 $D_0=\{a_{2i}:i\in\w\}$.

Let us show that for any distinct numbers $n,m$ the intersection
 $x_nA_{\ge n}\cap x_mA_{\ge m}$ is empty.
Otherwise there would exist two numbers $i\ge n$ and $j\ge m$ such that
 $x_na_{i}=x_ma_{j}$.
It follows from $x_n\ne x_m$ that $i\ne j$. We lose no generality
 assuming that $j>i$.
Then $x_na_{i}=x_ma_{j}$ implies that
$$a_{j}=x_m^{-1}x_na_{i}\in F_{j}^{-1}F_{j}A_{<j},$$ which contradicts
 the choice of $a_{j}$.

Let $\mathcal D\in\beta(X)$ be any free ultrafilter such that
 $D_0\in\mathcal D$ and $\mathcal D$
is not a P-point. To get such an ultrafilter, take $\mathcal D$ to be a
 cluster point of any
countable subset of $D^\pm_{0}\cap \beta(X)\setminus X$. Using the fact
 that $\mathcal D$ fails
 to be a P-point, we can take a decreasing sequence of sets
 $\{V_n:n\in\w\}\subset\mathcal D$
  having no pseudointersection in $\mathcal D$. The latter means that
 for every $D\in\mathcal D$
   the almost inclusion $D\subset^* V_n$ (which means that $D\setminus
 V_n$ is finite) holds only
   for finitely many numbers $n$.

For every $n\in\w$ let $D_n=V_n\cap A_{\ge n}\cap D_0$. We claim that
 the ultrafilter $\mathcal D$
 and the family $(D_n)_{n\in\w}$ satisfy the requirements of
 Definition~\ref{ultradetect}.

Take any subset $A\subset \w$ and consider the set $U_A=\bigcup_{n\in
 A}x_nD_n$.

First we verify that $U_A\cup yU_A\ne X$ for each $y\in X$. Find $m\in\w$ with $y^{-1}=x_m$ and take any odd number
$k>m$. We claim that $a_k\notin U_A\cup yU_A$. Otherwise, $a_{k}\in x_n D_n\cup x^{-1}_mx_nD_n$ for some $n\in A$. It
follows that $a_{k}=x_na_i$ or $a_{k}=x^{-1}_mx_na_i$ for some even $i\ge n$. If $k>i$, then both the equalities are
forbidden by the choice of $a_{k}\notin F^{-1}_kF_kA_{<k}\supset\{x_na_i,x^{-1}_mx_na_i\}$. If $k<i$, then those
equalities are formidden by the choice of
$$a_i\notin
F_i^{-1}F_iA_{<i}\supset\{x_n^{-1}a_k,x_n^{-1}x_m^{-1}a_k\}.$$Therefore, $U_A\cup yU_A\ne X$.

Next, given arbitrary $D\in\mathcal D$ we show that the set $S=\{n\in\w:x_nD\subset U_A\}$ is finite and lies in $A$.
First we show that $S\subset A$. Assuming the converse, we could find $n\in S\setminus A$. Then $x_n(D\cap D_n)\subset
x_nD\subset U_A=\bigcup_{m\in A}x_mD_m$, which is not possible because the set $x_nD_n$ misses the union $U_A$. Thus
$S\subset A$. Next, we show that $S$ is finite. By the choice of the sequence $(V_n)$, the set $F=\{n\in\w:D\cap
D_0\subset^* V_n\}$ is finite. We claim that $S\subset F$. Indeed, take any $m\in S$. It follows from $x_mD\subset
U_A=\bigcup_{n\in A}x_n D_n$ and $x_mA_{\ge m}\cap \bigcap_{n\ne m}x_nD_n=\emptyset$ that $$x_m(D\cap
D_0)\subset^*x_m(D\cap A_{\ge m})\subset x_mD_m\subset x_mV_m$$ and hence $m\in F$.
\end{proof}

\begin{theorem}\label{t3.3} Let $X$ be a group admitting a detecting
 ultrafilter $\mathcal D$.
For a maximal linked system $\A\in\lambda(X)$ the following conditions
 are equivalent:
\begin{enumerate}
\item the left shift $L_\A:G(X)\to G(X)$,  $L_\A:\F\mapsto \A\circ\F$,
 is continuous;
\item the left shift $l_\A:\lambda(X)\to\lambda(X)$, $l_\A:\mathcal
 L\mapsto \A\circ\mathcal L$, is continuous;
\item the left shift $l_\A:\lambda(X)\to\lambda(X)$ is continuous at
 the detecting ultrafilter $\mathcal D$;
\item $\A\in\lambda^\bullet(X)$.
\end{enumerate}
\end{theorem}

\begin{proof}
The implications $(1)\Ra(2)\Ra(3)$ are trivial while $(4)\Ra(1)$ follows from Theorem 7.1 \cite{G2} asserting that the
topological center of the semigroup $G(X)$ coincides with $G^\bullet(X)$. To prove that $(3)\Ra(4)$, assume that the
left shift $l_\A:\lambda(X)\to\lambda(X)$ is continuous at the detecting ultrafilter $\mathcal D$.

We need to show that $\A\in\lambda^\bullet(G)$. By Theorem~8.1 of \cite{G1}, it suffices to check that each set
$A\in\A$ contains a finite set $F\in\A$.

Since $\mathcal D$ is a detecting ultrafilter, there is a family of
 sets $\{D_x:x\in X\}\subset\mathcal D$ such that for every
 $D\in\mathcal D$
 the set $\{x\in X:xD\subset\bigcup_{a\in A}aD_A\}$ is finite and lies
 in $A$.

Consider the set $U_A=\bigcup_{x\in A}xD_x$ belonging to the product
 $\A\circ\mathcal D$.
 The continuity of the left shift $l_\A:\lambda(X)\to\lambda(X)$ at
 $\mathcal D$ yields us
 a set $D\in\mathcal D$, such that $l_\A(D^\pm)\subset U_A^\pm$. This
 means that $U_A\in\A\circ\mathcal L$
 for any maximal linked system $\mathcal L\in\lambda(X)$ that contains
 $D$.

The choice of $\mathcal D$ and $\{D_x\}_{x\in X}$ guarantees that
  $$S=\{x\in X:xD\subset U_A\}$$is
a finite subset lying in $A$. We claim that there is a maximal linked system $\tilde{\mathcal
 L}\in\lambda(X)$ such
that $D\in\tilde{\mathcal L}$ and $x^{-1}U_A\notin\tilde{\mathcal L}$
 for all $x\notin S$.
 Such a system $\tilde{\mathcal L}$ can be constructed as an
 enlargement of the linked system
 $$\mathcal L=\{D,X\setminus x^{-1}U_A:x\in X\setminus S\}.$$ The
 latter system is linked because of
 the definition of $S=\{x\in X:D\subset x^{-1}U_A\}$ and  the property
 (1) of the family $(D_x)_{x\in X}$
 from Definition~\ref{ultradetect}.

Take any maximal linked system $\tilde{\mathcal L}$ containing
 $\mathcal L$ and observe that $D\in\mathcal L$
 and $$\{x\in X:x^{-1}U_A\in\tilde{\mathcal L}\}=\{x\in
 X:x^{-1}U_A\in\mathcal L\}=S\subset A.$$
 Taking into account that $D\in\mathcal L$, we conclude that
 $\A\circ\mathcal L=l_\A(\mathcal L)\in U_A^\pm$
 and hence the set $S=\{x\in X:x^{-1}U_A\in\mathcal L\}\in\A$. This set
 $S$ is the required finite subset
 of $A$ belonging to $\A$.
\end{proof}

Combining Theorem~\ref{t3.3} with Lemma~\ref{l3.2} we obtain the main
 result of this section.

\begin{corollary}\label{t3.4} For any countable group $X$ the
 topological center of
the semigroup $\lambda(X)$ coincides with $\lambda^\bullet(X)$.
\end{corollary}

\begin{question} Is Theorem~\ref{t3.4} true for a group $X$ of
 arbitrary cardinality?
\end{question}

\section{The algebraic center of $\lambda(X)$}\label{center}

This section is devoted to studying the algebraic center of
 $\lambda(X)$.
We recall that the {\em algebraic center} of a semigroup $S$ consists
 of all
elements $s\in S$ that commute with all other elements of $S$. Such
 elements $s$ are called {\em central} in $S$.

\begin{lemma}\label{centralemma} Let $X$ be a group with the neutral
 element $e$. A maximal
 linked system $\A\in\lambda(X)$ is not central in $\lambda(X)$
 provided there are sets $S,T\subset X$ such that
\begin{enumerate}
\item $|T|=3$;
\item for each $A\in\A$ we get  $A\cap S\in\A$ and $|A\cap S|\ge2$;
\item there is a finite set $B\in\A$ such that $BS^{-1}\cap
 T^{-1}T\subset\{e\}$.
\end{enumerate}
\end{lemma}

\begin{proof} We claim that $\A$ does not commute with the maximal
 linked system
$\mathcal T=\{A\subset X:|A\cap T|\ge 2\}$. By (3), the maximal linked
 system $\A$ contains
a finite set $B\in\A$ such that $BS^{-1}\cap TT^{-1}\subset \{e\}$. By
 (2), we can assume that $B\subset S$
and $B$ is minimal in the sense that each $B'\subset B$ with $B'\in\A$
 is equal to $B$. By (2), $|B|\ge 2$.

Choose a family $\{T_b\}_{b\in B}$ of 2-element subsets of $T$ such
 that $\bigcup_{b\in B}T_b=T$.
Such a choice is possible because $|B|\ge2$.

 The union $\bigcup_{b\in B}bT_b$ belongs to $\A\circ\mathcal
 T=\mathcal T\circ\A$ and hence we
 can find a subset $D\in\mathcal T$ and a family $\{A_d\}_{d\in
 D}\subset\A$
 with $\bigcup_{d\in D}dA_d\subset \bigcup_{b\in B}bT_b$. By (2), we
 can assume that each
  $A_d\subset S$. Replacing $D$ by a smaller set, if necessary, we can
 assume that $D\subset T$
  and $|D|=2$. We claim that $A_d=B$ for all $d\in D$ and $T_b=D$ for
 all $b\in B$.

Indeed, take any $d\in D$ and any $a\in A_d$. Since $da\in \bigcup_{x\in D}xA_x\subset  \bigcup_{b\in B}bT_b$, there
are $b\in B$ and $t\in T_b$ with $da=bt$. Then $T^{-1}T\ni t^{-1}d=ba^{-1}\in BA_d^{-1}\subset BS^{-1}$. Taking into
account that $T^{-1}T\cap BS^{-1}\subset\{e\}$, we conclude that $t^{-1}d=ba^{-1}$ is the neutral element of $X$.
Consequently, $a=b\in B$ and $d=t\in T_b$. Since $a\in A_d$ was arbitrary, we get $\A\ni A_d\subset B$. The minimality
of $B\in\A$ implies that $A_d=B$. It follows from $d=t\in T_b$ for $d\in D$ that $D\subset T_b$. Since $|D|=|T_b|=2$,
we get $D=T_b$ for every $b\in B=A_d$. Consequently, $D=\bigcup_{b\in B}T_b=T$ which contradictions (1).
\end{proof}

By \cite[6.54]{HS}, for every group $X$ the algebraic center of
the
 semigroups $\beta(X)$
 coincides with the center of the group $X$. Consequently, the
 semigroup $\beta(X)\setminus X$
 contains no central elements. A similar result holds also for the
 semigroup $\lambda(X)$.

\begin{theorem}\label{algcent} For any countable infinite group $X$ the
 algebraic center of $\lambda(X)$
 coincides with the algebraic center of $X$.
\end{theorem}

\begin{proof} It is clear that all central elements of $X$ are central
 in $\lambda(X)$.
Now assume that a maximal linked system $\C\in\lambda(X)$ is a central
 element of the
semigroup $\lambda(X)$. Observe that the left shift
 $l_\C:\lambda(X)\to\lambda(X)$,
$l_\C:\mathcal X\mapsto \C\circ\mathcal X$ is continuous because it
 coincides with
the right shift $r_\C:\lambda(X)\to\lambda(X)$, $r_\C:\mathcal X\mapsto
 \mathcal X\circ\C$.
Consequently, $\C$ belongs to the topological center of $\lambda(X)$.
 Applying
Theorem~\ref{t3.4}, we conclude that $\C\in\lambda^\bullet(X)$. We
 claim that $\A$ is a principal ultrafilter.

Assuming the converse, consider the family $\C_0$ of minimal finite
 subsets in $\C$.
Since $\C\in\lambda^\bullet(X)$, the family $\C_0$ is finite and hence
 has finite union
$S=\cup\C_0$. Take any set $B\in\C_0$ and observe that $|B|\ge2$
 (because $\C$ is not a principal ultrafilter).

Since the group $X$ is infinite,  we can choose a 3-element subset
 $T\subset X$ such that
$T^{-1}T\cap BS^{-1}\subset\{e\}$. Now we see that the maximal linked
 system $\C$
satisfies the conditions of Lemma~\ref{centralemma} and hence is not
 central in $\lambda(X)$, which is a contradiction.
\end{proof}

We do not know if Theorem~\ref{algcent} is true for any infinite group
 $X$.

\begin{question} Let $X$ be an infinite group. Does the algebraic
 center of $\lambda(X)$
coincides with the algebraic center of $X$?
\end{question}

\begin{remark}\label{r4.4}
Theorem~\ref{algcent} certainly is not true for finite groups.
 According to \cite[\S \;6]{BGN}, for any group $X$ of cardinality $3\le|X|\le5$
 the semigroup $\lambda(X)$ contains a central element, which is not a
 principal ultrafilter.
\end{remark}

\begin{problem} Characterize (finite) abelian groups $X$ whose
 superextensions $\lambda(X)$
have central elements distinct from principal ultrafilters. Have all
 such groups $X$ cardinality $|X|\le 5$?
\end{problem}

It is interesting to remark that the semigroup $\lambda(X)$ contains
 many non-principal maximal
linked systems that commute with all ultrafilters.

\begin{proposition} Let $X$ be a group and $Y,Z\subset X$ be non-empty
 subsets such that $yz=zy$ for all
$y\in Y$, $z\in Z$. Then for any $\mathcal
 L\in\lambda^\bullet(Y)\subset\lambda^\bullet(X)$
and $\U\in\beta(Z)\subset\beta(X)$ we get $\mathcal
 L\circ\U=\U\circ\mathcal L$.
\end{proposition}

\begin{proof}

It is sufficient to prove that $\mathcal L\circ\U\subset\U\circ\mathcal L$. Let $\bigcup_{x\in L}x*U_x\in \mathcal
L\circ\U$. Without loss of generality we may assume that $L=\{x_1,\ldots,x_n\}$ is finite, $L\subset Y$ and
$U_{x_{i}}\subset Z$. Denote $V=U_{x_1}\cap\ldots\cap U_{x_n}\in\U$. Then
$$\bigcup_{x\in L}x*U_x=\bigcup_{x\in L}U_x*x\supset
V*L\in\U\circ\mathcal L.$$ It follows that $\bigcup_{x\in L}x*U_x\in\U\circ\mathcal L$ and the proof is complete.
\end{proof}


\begin{thebibliography}{10}
\bibitem[BGN]{BGN} T.~Banakh, V.~Gavrylkiv, O.~Nykyforchyn. {\em Algebra
 in superextensions of groups, I: zeros and commutativity}, preprint (arXiv:0802.1853).

\bibitem[BG$_3$]{BG3} T.~Banakh, V.~Gavrylkiv. {\em Algebra in
 superextension of groups, III: the minimal ideal}, preprint.


\bibitem[G$_1$]{G1} V.~Gavrylkiv. {\em The spaces of inclusion
 hyperspaces over noncompact
spaces} // Matem. Studii. {\bf 28}:1 (2007), 92--110.

\bibitem[G$_2$]{G2} V.~Gavrylkiv, {\em Right-topological semigroup
 operations on inclusion hyperspaces} // Matem. Studii (to appear).

\bibitem[H1]{Hind} N.~Hindman, {\em Finite sums from  sequences within
 cells of partition of $\IN$} // J.~Combin. Theory Ser. A {\bf 17}
 (1974), 1--11.

\bibitem[H2]{H2} N.~Hindman, {\em Ultrafilters and combinatorial number
 theory} // Lecture Notes in Math. {\bf 751} (1979), 49--184.

\bibitem[HS]{HS} N.~Hindman, D.~Strauss, {Algebra in the Stone-\v Cech
 compactification}, de Gruyter, Berlin, New York, 1998.

\bibitem[TZ]{TZ} A.~Teleiko, M.~Zarichnyi. Categorical Topology of
Compact Hausdofff Spaces, VNTL, Lviv, 1999.

\bibitem[P]{P} I.~Protasov. Combinatorics of Numbers, VNTL, Lviv, 1997.

\bibitem[vM]{vM} J.~van Mill,
Supercompactness and Wallman spaces, Math Centre Tracts. {\bf 85}. Amsterdam: Math. Centrum., 1977.


\end{thebibliography}
\end{document}